\documentclass[conf]{new-aiaa}

\title{Bayesian Algorithm for Collaborative Optimization with Application to Aircraft Design}

\author{
  Mohamed Ali Belhafnaoui
  \footnote{
    PhD Candidate,
    Department of Mathematics and Industrial Engineering,
    mohamed-ali.belhafnaoui@etud.polymtl.ca,
    AIAA Student Member.
  } and
  Youssef Diouane
  \footnote{
    Associate Professor,
    Department of Mathematics and Industrial Engineering,
    youssef.diouane@polymtl.ca,
    AIAA MDO TC Member.
  }
}
\affil{GERAD \& Polytechnique Montreal, Montreal, Quebec, H3T 1J4, Canada.}

\begin{document}
\maketitle
\begin{abstract}
  Collaborative Optimization (CO) is a multidisciplinary design optimization (MDO) framework that decomposes large-scale engineering problems into parallel, independently solvable subsystems coordinated by a system-level optimizer. Its practical utility is limited by the high frequency of expensive black-box disciplinary evaluations arising from the bi-level consistency constraints. This paper introduces BACO, a Bayesian Algorithm for Collaborative Optimization, which replaces the direct black-box calls at both levels with Gaussian process (GP) surrogates and acquisition function maximization. At the subsystem level, an acquisition function subject to GP-predicted feasibility constraints identifies the next evaluation point. At the system level, the same surrogate framework enforces consistency through predicted discrepancy constraints. This architecture reduces the number of true black-box evaluations required per major iteration. BACO is benchmarked against state-of-the-art CO variants on a Scalable MDO problem over 50 randomized instances. On this problem, BACO consistently achieves lower objective values and drives both constraint violation and interdisciplinary discrepancy to near-zero within the evaluation budget, outperforming all three CO variants across all tested DoE sizes. Further validation is conducted on a coupled aero-structural wing optimization problem based on the Common Research Model (CRM) geometry, where BACO identifies a feasible solution within 886 of 1000 allocated evaluations, recovering results physically consistent with active bending stress and tip deflection constraints. The BACO software, the state-of-the-art CO solvers, as well as standard MDO benchmarking problems are open-source and publicly available at \href{https://moebehfn.github.io/mdotoolbox/}{\textcolor{blue}{\texttt{moebehfn.github.io/mdotoolbox}}}.
\end{abstract}
\section{Nomenclature}
{
  \renewcommand\arraystretch{1.0}
  \noindent
  \begin{longtable*}{@{}l @{\quad=\quad} l@{}}
  CO   & Collaborative Optimization\\
   MDO  & Multidisciplinary Design Optimization\\
    BACO & Bayesian Algorithm for Collaborative Optimization\\
    ECO  & Enhanced Collaborative Optimization\\
    ICO  & Improved Collaborative Optimization\\
    MCO  & Modified Collaborative Optimization\\
    BO   & Bayesian Optimization\\
    CRM  & Common Research Model \\
    DoE  & Design of Experiments\\
    EI   & Expected Improvement\\
    GP   & Gaussian Process\\
    LHS  & Latin Hypercube Sampling\\
    SEGO & Super Efficient Global Optimization\\
    SoS  & System-of-Systems\\
    VoV  & Vector-of-Vectors\\

  \end{longtable*}
}
\setcounter{table}{0}
\newpage
\section{Introduction}\label{sec:intro}

\lettrine[realheight]{T}{he} design of complex aerospace systems requires coordination of multiple coupled physical disciplines, where design decisions in one discipline propagate through coupling variables to affect others. High-fidelity simulations of aerodynamics, structures, propulsion, and control systems are computationally expensive, often requiring hours or days per evaluation. This computational burden presents the challenge of how to systematically explore the design space and identify optimal configurations while minimizing the number of expensive function evaluations.

Multidisciplinary Design Optimization (MDO)~\cite{martins2013} addresses this challenge by integrating analysis and optimization across interacting subsystems. Collaborative Optimization (CO)~\cite{braun1996,braun1995} has emerged as particularly suited to large-scale problems, decomposing optimization into two hierarchical levels. A system-level coordinator that manages global objectives and enforces consistency among disciplines, and subsystem-level optimizations that solve local problems independently~\cite{martins2013}. However, its effectiveness is limited by the computational cost of subsystem optimizations, each requiring numerous evaluations of expensive black-box functions. Furthermore, standard CO implementations often struggle with load balancing. In a synchronous CO framework, the system-level coordinator must wait for all subsystems to terminate, even if some finish significantly earlier than others, leading to idle computational resources. This paper presents a survey of standard CO and its variants, including Modified CO (MCO)~\cite{alexandrov2002}, Enhanced CO (ECO)~\cite{roth2008}, and Improved CO (ICO)~\cite{li2026}, which often rely on gradient-based optimization.

Bayesian Optimization (BO)~\cite{jones1998,shahriari2016,frazier2018} has proven effective for expensive black-box optimization by constructing probabilistic surrogate models, often Gaussian Processes (GP)~\cite{rasmussen2005}, that guide the search with minimal function evaluations. An acquisition function balances exploration and exploitation~\cite{brochu2010}, enabling efficient convergence. BO has been successfully applied to multidisciplinary problems~\cite{priem2020,baars2024,difiore2023,cardoso2024}, demonstrating substantial computational cost reductions. In particular, the efficient global MDO approach~\cite{cardoso2024}  where  disciplinary surrogates are used to manage the coupling between non-linear solvers.

Despite these advances, the integration of BO within hierarchical MDO architectures like CO remains underexplored. Existing collaborative Bayesian optimization work~\cite{zhan2025,wang2025} addresses distributed optimization among peer agents, not the nested structure of CO where subsystem optimizations are subordinate to system-level decisions. This hierarchical coupling introduces unique challenges. Subsystem objectives depend on system-level targets that change across iterations, system-level evaluations require solving subsystem optimizations, introducing nested uncertainty, and consistency constraints~\cite{alexandrov2002} require specialized treatment within the BO framework.

As a primary contribution, we introduce the Bayesian Algorithm for Collaborative Optimization (BACO), which aims to explore a deep integration of BO into the CO framework to address computational inefficiency while preserving its modularity. This approach replaces traditional subsystem optimizers with independent GP surrogates that use acquisition functions to guide the selection of evaluation points for expensive black-box functions. The system-level coordinator maintains global consistency and manages the non-stationarity introduced by evolving targets to ensure that surrogate performance remains robust. This architecture leverages the exploration capabilities of BO to reduce the computational cost of hierarchical MDO without sacrificing the disciplinary autonomy inherent to the method.

The remainder of this paper is structured as follows. \Cref{sec:notations} defines the mathematical notation used throughout the document. \Cref{sec:background} reviews the CO framework and the BO methodology on which BACO is built. \Cref{sec:baco} presents the BACO formulation, including the dual-dataset discrepancy surrogate structure and the acquisition function problems at both levels. \Cref{sec:conv} defines the discrepancy and constraint violation metrics used to assess solution quality and termination. \Cref{sec:implementation} describes the implementation details and solver configurations used across all frameworks. \Cref{sec:bench-results} reports numerical results on a Scalable MDO benchmark problems, with statistical comparisons against the state-of-the-art. \Cref{sec:aeros-results} applies BACO to a coupled aero-structural wing optimization problem based on the Common Research Model (CRM)~\cite{vassberg2008} geometry. \Cref{sec:conclusion} summarizes the findings and identifies directions for future work.
\section{Notations}\label{sec:notations}

To maintain uniformity, this document employs a set of common mathematical notations. Vectors are represented as bold lower-case characters, such as~\(\bm{v}\). Matrices are represented as bold upper-case characters, such as~\(\bm{M}\). Digits printed in bold represent vectors of constant values, such as~\(\bm{1}={[1, 1, \ldots, 1]}^\top\).

Let~\(\bm{M} = [m_{ij}] \in \R^{m\times n}\) denote a matrix with entry~\(m_{ij}\) at row~\(i\), column~\(j\), having~\(i\)-th row~\(\bm{m}_{i,:}\) and~\(j\)-th column~\(\bm{m}_{:,j}\). A vector-of-vectors (VoV) of~\({\{v_{i}\}}_{1:l}\) is a vector~\(\bm{v}\in\Omega^n\) composed of sub-vectors~\(\bm{v}_i\in\Omega^{n_i}\), where~\(\bm{v}={[\bm{v}_1^\top, \bm{v}_2^\top, \ldots, \bm{v}_l^\top]}^\top\) and~\(\sum_{i=1}^l n_i = n\).

A scalar component of a VoV is~\({[\bm{v}_i]}_j\). For a given VoV,~\(\bm{v}\), the following set notations apply.
\begin{equation*}
  \bm{v}_{i\neq j} =
  \begin{cases}
    {\left[
        \bm{v}_2^\top,
        \bm{v}_3^\top,
        \cdots,
        \bm{v}_l^\top
    \right]}^\top & \text{for }\ j = 1\\
    {\left[
        \bm{v}_1^\top,
        \bm{v}_2^\top,
        \cdots,
        \bm{v}_{l-1}^\top
    \right]}^\top & \text{for }\ j = l\\
    {\left[
        \bm{v}_1^\top,
        \cdots,
        \bm{v}_{j-1}^\top,
        \bm{v}_{j+1}^\top,
        \cdots,
        \bm{v}_l^\top
    \right]}^\top & \text{for any }\ j\in\{2,\ldots,l-1\}
  \end{cases}
\end{equation*}
\begin{alignat*}{3}
  &\bm{v}_{i:j} &&= {\left[
      \bm{v}_i^\top,
      \bm{v}_{i+1}^\top,
      \cdots,
      \bm{v}_{j-1}^\top,
      \bm{v}_j^\top
  \right]}^\top\quad &&\text{for }\ 1\leq i<j \leq l\\
  &{\left[\bm{v}_k\right]}_{i:j} &&= {\left[
      {\left[\bm{v}_k\right]}_i,
      {\left[\bm{v}_k\right]}_{i+1},
      \cdots,
      {\left[\bm{v}_k\right]}_{j-1},
      {\left[\bm{v}_k\right]}_j
  \right]}^\top\quad &&\text{for }\ 1\leq i<j \leq n_k
\end{alignat*}

A scalar component of a regular vector retains the conventional notation~\(v_i\). The~\(n^\text{th}\) norm of a vector is formulated as~\(\|\bm{v}\|_n = {(\sum_{i=1}^n |v_i|^n)}^{1/n}\). For the second-order (Euclidean) norm, the subscript 2 is omitted, written as~\(\|\bm{v}\|\).

For optimization processes, the following conventions are used. Variables calculated at the system level are denoted with an overline (e.g.,~\(\overline{\bm{v}}\)), while those calculated at the subsystem level are denoted with an underline (e.g.,~\(\underline{\bm{v}}\)). The state of a variable,~\(\bm{v}\), at iteration instance~\(k\) is denoted as~\(\bm{v}^{(k)}\). An asterisk superscript,~\(\bm{v}^*\), marks the best solution found found.

Finally, function mapping is generalized. If~\(f: \R \to \R\), then~\(f(\bm{x} \in \R^n)\) implies element-wise application resulting in a vector in~\(\R^n\). Similarly, if~\(f: \R^n \to \R\), then~\(f(\bm{X} \in \R^{n \times m})\) implies application to each column of matrix~\(\bm{X}\), resulting in a vector in~\(\R^m\).

\section{Background and Context}\label{sec:background}

MDO frameworks are being developed to manage designing complex engineering systems. Among these, CO introduced a hierarchical decomposition strategy that enables distributed discipline-level optimization while enforcing system-level consistency constraints~\cite{martins2013}. In parallel, BO has emerged as a powerful methodology for optimizing expensive black-box functions by leveraging surrogate models to help identify the next expensive evaluation~\cite{frazier2018}. This section reviews the fundamental principles of CO and BO, providing the necessary context for integrating machine-learning-based methods into the CO framework.

\subsection{Collaborative Optimization}\label{subsec:collopt}

CO is a framework of multidisciplinary optimization that emerged as a potential solution for optimizing large-scale distributed analyses~\cite{braun1995}. This framework decomposes a problem hierarchically into a main system at the higher level and a set of subsystems (disciplines) at the lower level, coordinating between subcomponents to minimize the discrepancy (disagreement between designs) down to elimination. The use of this framework reduces information transfer between subcomponents and removes large iteration loops by optimizing the subcomponents independently and in parallel~\cite{martins2013, braun1995}.

Per~\cite{martins2013}, the system-level optimization problem~\eqref{prob:cosys} and the subsystem-level optimization problem~\eqref{prob:cosubsys} at subsystem~\(i\) are modeled as follows.
\begin{alignat}{3}
  &\min_{
    \zsys,
    \xsys,
    \ysys
  }\quad&&f(
    \zsys,
    \xsys,
    \ysys
  ) \tag{\(P_\text{sys}\)}\label{prob:cosys}\\
  &\text{subject to:}\quad&&\bm{c}(
    \zsys,
    \xsys,
    \ysys
  ) &&\geq \bm{0}\nonumber\\
  &  && J(
    \zsys,
    \xsys_i,
    \ysys,
    \zsub_i^*,
    \xsub_i^*
  ) &&= 0\quad \text{for all}\quad i\in\{1, \ldots, N\},\nonumber
\end{alignat}
and, for each subsystem~\(i\), the subsystem variables~\(\zsub_i^*\) and~\(\xsub_i^*\) are the solution of the following problem.
\begin{equation}
  \min_{
    \zsub_i,
    \xsub_i
  }\quad J\left(
    \zsys,
    \xsys_i,
    \ysys,
    \zsub_i,
    \xsub_i
  \right) \tag{\(P_i\)}\quad
  \text{subject to:}\quad\bm{g}_i\left(
    \ysys_{j\neq i},
    \zsub_i,
    \xsub_i
  \right)\geq\bm{0}.\label{prob:cosubsys}
\end{equation}

The discrepancy function~\(J\) ensures the design and coupling variables calculated at the system level match those from the subsystem level as follows.
\begin{equation}
  J(
    \zsys,
    \xsys_i,
    \ysys,
    \zsub_i,
    \xsub_i
  ) = J_i =
  \bigl\|\zsys-\zsub_i\bigr\|^2 +
  \bigl\|\xsys_i-\xsub_i\bigr\|^2 +
  \bigl\|\ysys_i-\bm{y}_i\bigl(
    \ysys_{j\neq i},
    \zsub_i,
    \xsub_i
  \bigr)\bigr\|^2.\label{eqn:discre}
\end{equation}

The components of the two levels of CO are as follows.

\begin{center}
  \begin{tabularx}{\textwidth}{@{}l @{\quad:\quad} X@{}}
    \(f\) & Objective function for problem~\eqref{prob:cosys} (e.g., fuel burn, profit, or emissions).\\

    \(\bm{c}\) & Constraint functions for problem~\eqref{prob:cosys} (e.g., total cost, weight, or endurance).\\

    \(\bm{g}_i\) & Constraint functions for problem~\eqref{prob:cosubsys} (e.g., glide ratio or specific fuel burn).\\

    \(\zsys\) (resp.,~\(\zsub_i\)) & Vector of the global design variables at the system level (resp., subsystem level).\\

    \(\xsys\) (resp.,~\(\xsub\)) & VoV of the local design variables at the system level (resp., subsystem level).\\

    \(\bm{y}\) & VoV of disciplinary analysis functions (or black-boxes). It is also assessed at the system level as a VoV, \(\ysys\), of inter-disciplinary coupling variables. This is usually the expensive part of the optimization.
  \end{tabularx}
\end{center}

\begin{figure}[!h]
  \centering
  \resizebox{.7\textwidth}{!}{%
    \begin{tikzpicture}[node distance=1cm]

  \def\layerdx{0.5cm}
  \def\layerdy{0.5cm}

  \begin{pgfonlayer}{back3}
    \begin{scope}[xshift=3*\layerdx, yshift=3*\layerdy]
      \node (inN)    [io                                      ] {\(\zsys, \xsys_N, \ysys_N\)};
      \node (ssn_yN) [bbeval,   below of=inN                  ] {\(\bm{y}_N\) Evaluation};
      \node (ssn_jN) [process,  below of=ssn_yN               ] {\(J_N, \bm{g}_N\) Evaluation};
      \node (dN)     [decision, below of=ssn_jN, yshift=-0.5cm] {Converged?};

      \draw [arrow     ] (inN)    -- (ssn_yN);
      \draw [arrow     ] (ssn_yN) -- (ssn_jN);
      \draw [arrow     ] (ssn_jN) -- (dN);
      \draw [arrow, red] (dN)     -| node[anchor=north] {no} ([xshift=-0.5cm]dN.west) |- (ssn_yN);
    \end{scope}
  \end{pgfonlayer}
  \begin{pgfonlayer}{back3card}
    \node [card, fit=(inN) (dN)] (cardN) {};
  \end{pgfonlayer}

  \begin{pgfonlayer}{back1}
    \begin{scope}[xshift=\layerdx, yshift=\layerdy]
      \node (in2)    [io                                      ] {\(\zsys, \xsys_2, \ysys_2\)};
      \node (ssn_y2) [bbeval,   below of=in2                  ] {\(\bm{y}_2\) Evaluation};
      \node (ssn_j2) [process,  below of=ssn_y2               ] {\(J_2, \bm{g}_2\) Evaluation};
      \node (d2)     [decision, below of=ssn_j2, yshift=-0.5cm] {Converged?};

      \draw [arrow     ] (in2)    -- (ssn_y2);
      \draw [arrow     ] (ssn_y2) -- (ssn_j2);
      \draw [arrow     ] (ssn_j2) -- (d2);
      \draw [arrow, red] (d2)     -| node[anchor=north] {no} ([xshift=-0.5cm]d2.west) |- (ssn_y2);
    \end{scope}
  \end{pgfonlayer}
  \begin{pgfonlayer}{back1card}
    \node [card, fit=(in2) (d2)] (card2) {};
  \end{pgfonlayer}

  \node at ($ (card2.north west)!0.125!(cardN.north west) $)
  [anchor=west] {\Large\(\bm{\cdots}\)};

  \node (in1)    [io                                      ] {\(\zsys, \xsys_1, \ysys_1\)};
  \node (ssn_y1) [bbeval,   below of=in1                  ] {\(\bm{y}_1\) Evaluation};
  \node (ssn_j1) [process,  below of=ssn_y1               ] {\(J_1, \bm{g}_1\) Evaluation};
  \node (d1)     [decision, below of=ssn_j1, yshift=-0.5cm] {Converged?};

  \draw [arrow     ] (in1)    -- (ssn_y1);
  \draw [arrow     ] (ssn_y1) -- (ssn_j1);
  \draw [arrow     ] (ssn_j1) -- (d1);
  \draw [arrow, red] (d1)     -| node[anchor=north] {no} ([xshift=-0.5cm]d1.west) |- (ssn_y1);

  \begin{pgfonlayer}{maincard}
    \node [card, fit=(in1) (d1)] (card1) {};
  \end{pgfonlayer}


  \node (ins) [io, right of=inN, xshift=4.5cm] {\(\zsub, \xsub\)};

  \node (s_y)  [bbeval,   below of=ins                                    ] {\(\{\bm{y}_i\}_{1:N}\) Evaluation};
  \node (s_j)  [process,  below of=s_y                                    ] {\(f, \bm{c}, \Jtot\) Evaluation};
  \node (s_d)  [decision, below of=s_j, yshift=-0.5cm                     ] {Converged?};
  \node (exit) [decision, below of=s_d, minimum height=0.75cm, yshift=-2.25cm] {Stop?};
  \node (outp) [io, below of=exit, yshift=-0.5cm                          ] {Optimization Output};

  \draw [arrow     ] (ins) -- (s_y);
  \draw [arrow     ] (s_y) -- (s_j);
  \draw [arrow     ] (s_j) -- (s_d);
  \draw [arrow     ] (exit) -- node[anchor=east] {yes} (outp);

  \draw [arrow, red] (s_d)       -| node[anchor=west]       {no}  ([xshift=0.5cm]s_d.east) |- (s_y);
  \draw [arrow     ] (s_d.south) -| node[anchor=north west] {yes} (exit.north);

  \draw [arrow] (d1.south)
  -- ([yshift=-.5cm]d1.south) node[anchor=north] {yes} -| ($ (inN.east)!0.5!(ins.west) $) -- (ins.west);
  \draw [arrow, red] (exit)
  -| ([yshift=-1.5cm,xshift=-1cm]d1.west) |- (in1.west);

  \begin{pgfonlayer}{back1}
    \draw [arrow] (d2.south)
    -- ([yshift=-.5cm]d2.south) node[anchor=north] {yes} -|
    node[anchor=south, rotate=90, xshift=\layerdy] {\(\bm{\cdots}\)}
    ($ (inN.east)!0.5!(ins.west) $)
    -- (ins.west);

    \draw [arrow, red] (exit)
    -| ([yshift=-1.5cm,xshift=-1cm]d1.west) |- (in2.west);
  \end{pgfonlayer}

  \begin{pgfonlayer}{back3}
    \draw [arrow] (dN.south)
    -- ([yshift=-.5cm]dN.south) node[anchor=north] {yes} -| ($ (inN.east)!0.5!(ins.west) $) -- (ins.west);

    \draw [arrow, red] (exit) -| node[anchor=north] {no} ([yshift=-1.5cm,xshift=-1cm]d1.west)
    |- node[anchor=south, rotate=-90, xshift=\layerdy] {\(\bm{\cdots}\)} (inN.west);
  \end{pgfonlayer}

\end{tikzpicture}
  }%
  \caption{Collaborative Optimization Data Flow}\label{fig:codiag}
\end{figure}
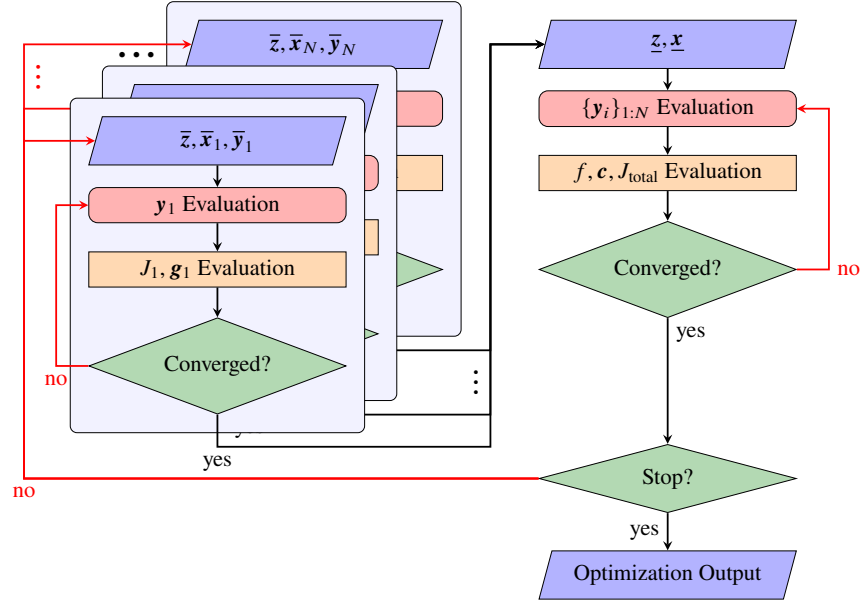

The structural decomposition and the iterative information exchange between the system and subsystem levels are illustrated in~\Cref{fig:codiag}.

A detailed description of the CO framework, as inspired by~\cite[Algorithm 3]{martins2013}, is given in~\Cref{alg:co}.

\begin{algorithm}[!h]
  \caption{A Standard Collaborative Optimization Framework}\label{alg:co}
  \DontPrintSemicolon
  \SetKwInOut{Input}{Inputs}\SetKwInOut{Output}{Outputs}

  \Input{
    \ul{Functions}:\\{
      \qquad
      \makebox[3.5cm][l]{\(f\): Objective function}
      \makebox[3.5cm][l]{\(\bm{c}\): Vector of system constraints}\\
      \qquad
      \makebox[3.5cm][l]{\(\bm{y}\): VoV of \({\{\bm{y}_i\}}_{1:N}\)}
      \makebox[3.5cm][l]{\(\bm{g}\): VoV of \({\{\bm{g}_i\}}_{1:N}\)}\\
    }
    \ul{Initial decision variables}:\\{
      \qquad
      \makebox[3.5cm][l]{\(\zsys^{(0)}\): Vector}
      \makebox[3.5cm][l]{\(\xsys^{(0)}\): VoV of \({\{\xsys_i\}}_{1:N}\)}
      \makebox[3.5cm][l]{\(\ysys^{(0)}\): VoV of \({\{\ysys_i\}}_{1:N}\)}\\
      \qquad
      \makebox[4.6cm][l]{\(\zsub^{(0)}\): Vector of identical size to \(\zsys\)}
      \makebox[4.6cm][l]{\(\xsub^{(0)}\): Vector of identical size to \(\xsys\)},\\
      \qquad\(\bm{\epsilon}\): Vector of thresholds.
    }
  }

  \(k \gets 0\)

  \Repeat{a stopping criterion is met}{
    \ForEach{subsystem \(i \in \left\{ 1, \ldots, N \right\}\)}{

      \textcolor{green}{\tcp{--- Subsystem Level, Parallel Execution ---}}

      Solve the problem~\eqref{prob:cosubsys}:\vspace{-.25cm}\[
        {(
            \zsub_i^{(k+1)},
            \xsub_i^{(k+1)}
        )}
        \in
        \argmin_{
          \zsub_i,
          \xsub_i
        }
        \left\{
          J{(
              \zsys^{(k)},
              \xsys_i^{(k)},
              \ysys^{(k)},
              \zsub_i,
              \xsub_i
          )}
          :
          \bm{g}_i{(
              \ysys_{j\neq i}^{(k)},
              \zsub_i,
              \xsub_i
          )}\geq\bm{0}
        \right\}
      \]\vspace{-.5cm}

    }

    \textcolor{green}{\tcp{--- System Level ---}}

    Solve the problem~\eqref{prob:cosys}:\vspace{-.25cm}\[
      {(
          \zsys^{(k+1)},
          \xsys^{(k+1)},
          \ysys^{(k+1)}
      )}
      \in
      \argmin_{
        \zsys,
        \xsys,
        \ysys
      }
      \left\{
        f{(
            \zsys,
            \xsys,
            \ysys
        )}
        :~
        \bm{c}{(
            \zsys,
            \xsys,
            \ysys
        )}\geq\bm{0}
        ~\text{and}~
        J{(
            \zsys,
            \xsys_i,
            \ysys,
            \zsub_i^{(k+1)},
            \xsub_i^{(k+1)}
        )} = 0\ \text{for}\ i = 1:N
      \right\}
    \]\vspace{-.5cm}

    \(k \gets k + 1\)
  }
\end{algorithm}

As depicted, the coordination follows a cyclic bi-level process. The system optimizer problem~\eqref{prob:cosys} dictates the target design and coupling variables~\((\zsys, \xsys, \ysys)\) to the lower level. In response, each subsystem operates independently in parallel to identify local solutions~\((\zsub_i^*, \xsub_i^*)\) that satisfy local constraints~\(\bm{g}_i\) while minimizing the discrepancy~\(J\) relative to the system targets. These local optima are returned to the system level, which then updates the global targets to minimize~\(f\) and satisfy~\(\bm{c}\). The specific procedural logic for this exchange is detailed in~\Cref{alg:co}.

CO was developed in collaboration with The National Aeronautics and Space Administration (NASA)~\cite{braun1995} and has since been used in many aerospace engineering applications, such as aircraft design~\cite{sobieski1996} and launch vehicle design~\cite{braun1996}. These applications require expensive black-boxes, such as Computational Fluid Dynamics, in the optimization process. While CO effectively decomposes complex multidisciplinary design problems, its performance remains tied to frequent evaluations of expensive black-boxes. This motivates methods that minimize the number of function calls. Bayesian optimization~\cite{frazier2018} addresses this by using surrogates to guide sampling, thereby reducing reliance on direct evaluations. Hybrid approaches within each CO iteration can further mitigate computational costs, leveraging surrogate-based search to maintain coordination benefits while substantially limiting the frequency of expensive solver calls.

CO's decomposition into parallel subsystem optimizations offers potential computational advantages for distributed computing environments~\cite{braun1995,martins2013}. However, practical implementation faces challenges including communication overhead from transmitting design targets~\((\zsys, \xsys, \ysys)\) to subsystems and returning optimized solutions~\((\zsub_i^*, \xsub_i^*)\) at each iteration. For problems with high-dimensional coupling variables or many subsystems, this data transfer can dominate wall-clock time, particularly for geographically distributed computational resources. Load balancing presents additional complications where subsystem optimization times may vary significantly depending on problem difficulty and active constraint sets, causing idle time for faster subsystems awaiting slower ones.

The computational burden of repeated black-box evaluations within each subsystem optimization has motivated surrogate-based extensions to CO. Ref.~\cite{baars2024} proposed replacing each discipline with an independent GP surrogate and selecting enrichment points via Thompson sampling~\cite{baars2024}, avoiding the need for a full multidisciplinary analysis at each iteration. This partitioned surrogate approach demonstrates that discipline-level surrogates can reduce evaluation counts in strongly coupled systems, but it relies on Thompson sampling as the infill criterion, which does not produce an explicit acquisition function amenable to gradient-based maximization. The present work builds on this direction by retaining the discipline-level surrogate structure while replacing the infill criterion with an acquisition function framework compatible with the CO consistency constraints.

Ref.~\cite{alexandrov2002} proved that the standard CO formulation is mathematically ill-posed because discrepancy constraints (such as $J = 0$) create a singularity that causes gradient-based optimizers to fail, an issue inherent to the degeneracy of its bi-level structure~\cite{luo1996,dempe2002,kim2020}. To resolve this, researchers have developed penalty-based and relaxation-based approaches to smooth or penalize these equilibrium constraints~\cite{alexandrov2002,roth2008,li2026}. MCO~\cite{alexandrov2002} reformulates the system level as an unconstrained problem, eliminating the singular Jacobian, and restructures the subsystem problems so that the multipliers do not converge to zero. MCO retains a weakness inherent to the bi-level structure: the subsystem responses remain non-smooth functions of the system-level targets due to changes in the active set of subsystem constraints.

Ref.~\cite{roth2008} introduces ECO, which shares linear models of the subsystem constraints with all other disciplinary design teams, communicating the source of non-smoothness and preventing conflicting struggles over shared variables. The system level becomes an unconstrained minimization problem focused on compatibility between subsystems, while subsystems incorporate a quadratic model of the global objective and thereby gain enhanced design authority.

Ref.~\cite{li2026} proposes ICO which incorporates contemporary optimization techniques combined with ensemble surrogate models~\cite{chen2024} to approximate implicit relationships. ICO with dynamic penalty (ICO-DP) transforms the system-level equality constraints into gradually convergent inequality constraints using a hypersphere relaxation mechanism, with a dynamic penalty function adjusted iteratively to guide the solution toward the feasible region. The approach is limited in high-dimensional problems, as Kriging surrogate construction scales poorly with dimension and per-iteration global optimization for the penalty factor increases overall cost.

\subsection{Bayesian Optimization Framework}\label{subsec:bayesopt}
This section covers optimization problems in the following form
\begin{equation}
  \min_{\bm{x}\in X}~
  f(\bm{x})\quad
  \text{subject to:}\quad
  \bm{x}\in\Omega=\{\bm{x}:\bm{h}(\bm{x})=\bm{0},~\bm{c}(\bm{x})\ge\bm{0}\},
  \tag{\(P\)}\label{prob:basic}
\end{equation}
where~\(\bm{x}={[x_1, x_2, \ldots, x_n]}^\top\) is a vector that represents the decision variables,~\(f:\mathcal{X}\rightarrow\R\) is the objective function,~\(\mathcal{X}\) is the domain,~\(\Omega\) is the set of points allowed as a solution, and~\(\bm{h}:\mathcal{X}\rightarrow\R^m\) and~\(\bm{c}:\mathcal{X}\rightarrow\R^k\) are the vectors of constraints.

BO~\cite{frazier2018,jones1998,shahriari2016} is a sequential, GP-based framework for global optimization of expensive black-box problems. A GP~\cite{rasmussen2005} prior is commonly employed to capture uncertainty in the objective function, enabling principled trade-offs between exploration and exploitation~\cite{cordelier2025}.

The GP regression framework is presented in generic terms using~\(s:\mathcal{X}\to\R\) to denote any scalar function requiring surrogate modeling. In the optimization context,~\(s(\bm{x})\) may represent an objective function~\(f(\bm{x})\), a constraint function~\(c_j(\bm{x})\), or any other expensive-to-evaluate black-box quantity. Subsequent sections apply this framework to specific function types.
Rather than fitting a single parametric model, a GP~\cite{rasmussen2005} places a probability distribution directly over functions. A GP is fully specified by its mean function~\(m(\bm{x})\) and covariance (kernel) function~\(k(\bm{x}, \bm{x}')\)~\cite{rasmussen2005}. Without loss of generality, assume~\(m(\bm{x}) = 0\). The GP prior is then
\begin{equation*}
  s(\bm{x}) \sim \mathcal{GP}\left(0, k(\bm{x}, \bm{x}')\right).
\end{equation*}
Once observations~\(\D=\{\bm{X}, \bm{s}\}\) are available, also know as Design of Experiments (DoE), the GP prior is updated to a posterior distribution~\(\mathcal{N}(\mu_s(\bm{x}), \sigma_s^2(\bm{x}))\). The posterior mean~\(\mu_s(\bm{x})\) serves as a weighted sum of the training observations, while the variance~\(\sigma_s^2(\bm{x})\) quantifies the epistemic uncertainty. At observed points, the model becomes certain, while at unobserved points, it retains calibrated uncertainty. This update follows the rules of conditional Gaussian distributions where the posterior mean and variance are analytically computable~\cite{rasmussen2005}.

At each iteration, BO builds a GP model using all available evaluations (current DoE), and an acquisition function~\cite{jones1998,shahriari2016} defined over the surrogate is maximized to select the next evaluation point. This decouples the expensive function evaluation from the search for a promising candidate, allowing near-optimal solutions to be located with a small number of evaluations~\cite{jones1998,priem2020}. The acquisition function~\(\alpha: \mathcal{X} \to \R\)~\cite{jones1998,frazier2018,shahriari2016} guides the optimization by balancing exploration and exploitation. Let~\(s^* = \min_{i = 1,\ldots,N} s(\bm{x}_i)\) denote the best observed value. A popular acquisition function choice is the Expected Improvement (EI) criterion~\cite{jones1998}, which  assesses the average amount by which the result would improve upon the current best~\(s^*\) if~\(s\) is evaluated at~\(\bm{x}\)
\begin{equation}\label{eqn:ei}
  \alpha_\text{EI}(\bm{x}) = (s^* - \mu_s(\bm{x})) \Phi(Z(\bm{x})) + \sigma_s(\bm{x}) \phi(Z(\bm{x})),
\end{equation}
where~\(Z(\bm{x}) = (s^* - \mu_s(\bm{x}))/\sigma_s(\bm{x})\), and~\(\Phi\) and~\(\phi\) denote the standard normal cumulative distribution and probability density functions.
To address numerical vanishing of gradients in regions far from the optimum, the Log Expected Improvement (LogEI)~\cite{ament2023} is used
\begin{equation}\label{eqn:logei}
  \alpha_\text{LogEI}(\bm{x}) = \texttt{log\_h}(Z(\bm{x})) + \log\sigma_s(\bm{x}).
\end{equation}

Because the logarithm is a monotone increasing function, LogEI and EI share the same maximizer: the point~\(\bm{x}^*\) that maximizes EI also maximizes LogEI. The difference is purely numerical: by working in log space, gradients that would be zero in the original scale become finite and informative, making the acquisition maximization well-conditioned even in difficult regions. The numerically stable helper function~\(\texttt{log\_h}: \mathbb{R} \to \mathbb{R}\)~\cite{ament2023} handles three regimes of~\(Z(\bm{x})\) to maintain precision throughout:
\begin{equation}\label{eqn:log_h}
  \texttt{log\_h}(z) =
  \begin{cases}
    \log\left(\phi(z) + z\Phi(z)\right) & z > -1 \\
    -z^2/2 - c + \log\left(1 - \texttt{log1mexp}\left(\log\left(\texttt{erfcx}\left(\dfrac{-z}{\sqrt{2}}\right)|z|\right) + c\right)\right) & -1/\sqrt{\varepsilon} < z \le -1 \\
    -z^2/2 - c - 2\log|z| & z \le -1/\sqrt{\varepsilon}
  \end{cases},
\end{equation}
with~\(c = \tfrac{1}{2}\log(2\pi)\),~\(\varepsilon > 0\) is the machine precision,~\(\texttt{log1mexp}(z) = \log(1 - \exp(z))\), and~\(\texttt{erfcx}(z) = \exp(z^2)\operatorname{erfc}(z)\) the scaled complementary error function~\cite{ament2023}. LogEI is the preferred formulation whenever gradient-based acquisition maximization is used~\cite{ament2023,cordelier2025}.

Many engineering optimization problems involve constraints that are themselves expensive to evaluate or available only as black-box outputs.
Super Efficient Global Optimization (SEGO)~\cite{sasena2002} is a BO framework for expensive black-box constrained optimization problems of the form~\eqref{prob:basic}. \Cref{alg:sego} presents the SEGO procedure following~\cite{sasena2002,cordelier2025}.

\begin{algorithm}[!h]
  \caption{SEGO: Super Efficient Global Optimization}\label{alg:sego}
  \DontPrintSemicolon
  \SetKwInOut{Input}{Inputs}\SetKwInOut{Output}{Outputs}\SetKwInOut{Initialize}{Initialize}

  \Input{%
    Initial DoE: \(
      \D^{(0)} = \{
        \bm{X}^{(0)},
        f(\bm{X}^{(0)}),
        \bm{c}(\bm{X}^{(0)}),
        \bm{h}(\bm{X}^{(0)})
      \}
    \) \textcolor{green}{\tcp*{Multifunction DoE}}
  }

  \Output{Optimum solution: \(f^*, \bm{x}^*\)}

  \(k \gets 0\);

  \Repeat{a stopping criterion is satisfied}{
    Fit GP surrogates for \(f, \bm{c}, \bm{h}\) using \(\D^{(k)}\).

    Solve acquisition problem:\vspace{-.25cm}\[
      \bm{x}^{(k+1)} \in \argmax_{\bm{x}} \{
        \alpha_f(\bm{x})
        :
        \bm{\mu}_c(\bm{x}) \geq \bm{0}
        ,
        \bm{\mu}_h(\bm{x}) = \bm{0}
      \}
    \]\vspace{-.5cm}

    Evaluate true functions and update DoE:\vspace{-.25cm}\[
      \D^{(k+1)} = \D^{(k)} \cup \{
        \bm{x}^{(k+1)},
        f(\bm{x}^{(k+1)}),
        \bm{c}(\bm{x}^{(k+1)}),
        \bm{h}(\bm{x}^{(k+1)})
      \}
    \]\vspace{-.5cm}

    Update best observed solution:\vspace{-.25cm}\[
      \bm{x}^* \in \argmin_{\bm{x} \in \bm{X}^{(k+1)}}
      \{
        f(\bm{x}):
        \bm{c}(\bm{x})\ge\bm{0},
        \bm{h}(\bm{x})=\bm{0}
      \};
      \quad
      f^* \gets f(\bm{x}^*)
    \]\vspace{-.5cm}

    \(k \gets k+1\);
  }
\end{algorithm}

In~\Cref{alg:sego}, when the feasible region is unknown at initialization, space-filling designs such as Latin Hypercube Sampling (LHS)~\cite{mckay1979} provide initial coverage. The Upper Trust Bound~\cite{priem2020b} approach addresses this by constructing a unified criterion for mixed constraint types, including hidden constraints~\cite{ledigabel2024} that may cause simulation failure without returning a function value.


\section{BACO: A Bayesian Algorithm for Collaborative Optimization}\label{sec:baco}

The system level problem in CO is defined as a function of~\(\zsys\),~\(\xsys\), and~\(\ysys\), with~\(\zsub_i\) and~\(\xsub_i\) being given constants fed from the subsystem problems. Similarly, the subsystem problem is only defined as a function of~\(\zsub_i\) and~\(\xsub_i\), with~\(\zsys\),~\(\xsys_i\), and~\(\ysys\) being given constants fed from the system problem. While building the GP surrogates, it is crucial to include these constants as inputs for the discrepancy functions~\(J_i\), as these vary with each iteration. Furthermore, the GP surrogate of the discrepancy~\(J_i\) can be enriched at both the subsystem and system levels within each system-level iteration, since~\(J_i\) appears both as the objective of the \(i\)-th subsystem problem and as a constraint of the system level problem.

Given this, the system level problem~\eqref{prob:cosys} and the subsystem level problem~\eqref{prob:cosubsys} become the following.
\begin{alignat}{5}
  &\max_{
    \zsys,
    \xsys,
    \ysys
  }\quad && \alpha_f\left(
    \zsys,
    \xsys,
    \ysys
  \right)
  \tag{\(\widetilde{P}_\text{sys}\)}\label{prob:bacosys}\\
  &\text{subject to:}\quad
  &&\bm{\mu_c}\left(
    \zsys,
    \xsys,
    \ysys
  \right) && \geq \bm{0} \nonumber\\
  &   &&\mu_{J_i}\left(
    \zsys,
    \xsys_i,
    \ysys,
    \zsub_i^*,
    \xsub_i^*
  \right) &&= 0\qquad \text{for all}\quad i\in\{1\cdots N\},\nonumber
\end{alignat}
and, following problems~(\ref{prob:cosys},~\ref{prob:cosubsys}), for each subsystem~\(i\), the subsystem variables~\(\zsub_i^*\) and~\(\xsub_i^*\) are the solution of the following problem.
\begin{alignat}{5}
  \max_{\zsub_i, \xsub_i}
  \alpha_{J_i}\left(
    \zsys,
    \xsys_i,
    \ysys,
    \zsub_i,
    \xsub_i
  \right)
  \quad\text{subject to:}\quad
  \bm{\mu}_{\bm{g}_i}\left(
    \ysys_{j\neq i},
    \zsub_i,
    \xsub_i
  \right) &&\geq\bm{0}.\tag{\(\widetilde{P}_i\)}\label{prob:bacosubsys}
\end{alignat}

This formulation induces three families of datasets. The system level functions are recorded as follows.
\begin{equation*}
  \D_\text{sys} = \left\{
    \overline{\bm{Z}},
    \overline{\bm{X}},
    \overline{\bm{Y}},
    f\left(
      \overline{\bm{Z}},
      \overline{\bm{X}},
      \overline{\bm{Y}}
    \right),
    \bm{c}\left(
      \overline{\bm{Z}},
      \overline{\bm{X}},
      \overline{\bm{Y}}
    \right)
  \right\}.
\end{equation*}

Two separate datasets are maintained for each discrepancy function~\(J_i\), one populated by system-level evaluations and one by subsystem-level evaluations.
\begin{equation*}
  \overline{\D_{J_i}} = \left\{
    \overline{\bm{Z}},
    \overline{\bm{X}}_i,
    \overline{\bm{Y}},
    \underline{\bm{Z}}_i,
    \underline{\bm{X}}_i,
    J_i\left(
      \overline{\bm{Z}},
      \overline{\bm{X}}_i,
      \overline{\bm{Y}},
      \underline{\bm{Z}}_i,
      \underline{\bm{X}}_i
    \right)
  \right\}
  \quad,\quad
  \underline{\D_{J_i}} = \left\{
    \overline{\bm{Z}},
    \overline{\bm{X}}_i,
    \overline{\bm{Y}},
    \underline{\bm{Z}}_i,
    \underline{\bm{X}}_i,
    J_i\left(
      \overline{\bm{Z}},
      \overline{\bm{X}}_i,
      \overline{\bm{Y}},
      \underline{\bm{Z}}_i,
      \underline{\bm{X}}_i
    \right)
  \right\}.
\end{equation*}

\begin{figure}[!h]
  \centering
  \resizebox{.65\textwidth}{!}{%
    \begin{tikzpicture}[node distance=1.5cm]

  \def\layerdx{0.5cm}
  \def\layerdy{0.5cm}

  \begin{pgfonlayer}{back3}
    \begin{scope}[xshift=3*\layerdx, yshift=3*\layerdy]
      \node (inN_b)    [io,      minimum width=4cm]                                  {\(\underline{\D}_{J_N}, \underline{\D}_{N}\)};
      \node (gpN_b)    [process, minimum width=4cm, below of=inN_b,  yshift=.5cm]    {Fitting GPs of~\(J_N, \bm{g}_N\)};
      \node (pN_b)     [process, minimum width=4cm, below of=gpN_b,  yshift=.5cm]    {Solving~\(\widetilde{P}_N\)};
      \node (yN_b)     [bbeval,  minimum width=4cm, below of=pN_b,   yshift=.5cm]    {\(\bm{y}_N\)~Evaluation};
      \node (jN_b)     [process, minimum width=4cm, below of=yN_b,   yshift=.5cm]    {\(\bm{g}_N\)~Evaluation};
      \node (doeN_b)   [process, minimum width=4cm, below of=jN_b,   yshift=.5cm]    {Update~\(\underline{\D}_{J_N}, \underline{\D}_{N}\)};
      \draw [arrow] (inN_b)  -- (gpN_b);
      \draw [arrow] (gpN_b)  -- (pN_b);
      \draw [arrow] (pN_b)   -- (yN_b);
      \draw [arrow] (yN_b)   -- (jN_b);
      \draw [arrow] (jN_b)   -- (doeN_b);
    \end{scope}
  \end{pgfonlayer}
  \begin{pgfonlayer}{back3card}
    \node [card, fit=(inN_b) (doeN_b)] (cardN) {};
  \end{pgfonlayer}

  \begin{pgfonlayer}{back1}
    \begin{scope}[xshift=\layerdx, yshift=\layerdy]
      \node (in2_b)    [io,      minimum width=4cm]                                  {\(\underline{\D}_{J_2}, \underline{\D}_{2}\)};
      \node (gp2_b)    [process, minimum width=4cm, below of=in2_b, yshift=.5cm]     {Fitting GPs of~\(J_2, \bm{g}_2\)};
      \node (p2_b)     [process, minimum width=4cm, below of=gp2_b, yshift=.5cm]     {Solving~\(\widetilde{P}_2\)};
      \node (y2_b)     [bbeval,  minimum width=4cm, below of=p2_b,  yshift=.5cm]     {\(\bm{y}_2\)~Evaluation};
      \node (j2_b)     [process, minimum width=4cm, below of=y2_b,  yshift=.5cm]     {\(\bm{g}_2\)~Evaluation};
      \node (doe2_b)   [process, minimum width=4cm, below of=j2_b,  yshift=.5cm]     {Update~\(\underline{\D}_{J_2}, \underline{\D}_{2}\)};
      \draw [arrow] (in2_b)  -- (gp2_b);
      \draw [arrow] (gp2_b)  -- (p2_b);
      \draw [arrow] (p2_b)   -- (y2_b);
      \draw [arrow] (y2_b)   -- (j2_b);
      \draw [arrow] (j2_b)   -- (doe2_b);
    \end{scope}
  \end{pgfonlayer}
  \begin{pgfonlayer}{back1card}
    \node [card, fit=(in2_b) (doe2_b)] (card2) {};
  \end{pgfonlayer}

  \node at ($ (card2.north west)!0.125!(cardN.north west) $)
  [anchor=west] {\Large\(\bm{\cdots}\)};

  \node (in1)      [io,      minimum width=4cm]                                    {\(\underline{\D}_{J_1}, \underline{\D}_{1}\)};
  \node (ss1_gp)   [process, minimum width=4cm, below of=in1,     yshift=.5cm]     {Fitting GPs of~\(J_1, \bm{g}_1\)};
  \node (ss1_p1)   [process, minimum width=4cm, below of=ss1_gp,  yshift=.5cm]     {Solving~\(\widetilde{P}_1\)};
  \node (ss1_y1)   [bbeval,  minimum width=4cm, below of=ss1_p1,  yshift=.5cm]     {\(\bm{y}_1\)~Evaluation};
  \node (ss1_j1)   [process, minimum width=4cm, below of=ss1_y1,  yshift=.5cm]     {\(\bm{g}_1\)~Evaluation};
  \node (ss1_doe1) [process, minimum width=4cm, below of=ss1_j1,  yshift=.5cm]     {Update~\(\underline{\D}_{J_1}, \underline{\D}_{1}\)};

  \draw [arrow] (in1)     -- (ss1_gp);
  \draw [arrow] (ss1_gp)  -- (ss1_p1);
  \draw [arrow] (ss1_p1)  -- (ss1_y1);
  \draw [arrow] (ss1_y1)  -- (ss1_j1);
  \draw [arrow] (ss1_j1)  -- (ss1_doe1.north);

  \begin{pgfonlayer}{maincard}
    \node [card, fit=(in1) (ss1_doe1)] (card1) {};
  \end{pgfonlayer}

  \node (ins) [io, right of=inN_b, xshift=4.5cm, trapezium stretches body=true] {
    \(\overline{\D}_{J_1}, \overline{\D}_\text{sys}\)
  };

  \node (s_gp) [process, minimum width=4cm, below of=ins,  yshift=.5cm] {Fitting GPs of~\(\{J_i\}_{1:N}, f, \bm{c}\)};
  \node (ps)   [process, minimum width=4cm, below of=s_gp, yshift=.5cm] {Solving~\(\widetilde{P}_\text{sys}\)};

  \draw [arrow] (ins)  -- (s_gp);
  \draw [arrow] (s_gp) -- (ps);

  \node (s_y)   [bbeval,  minimum width=4cm, below of=ps,   yshift=.5cm] {\(\{\bm{y}_i\}_{1:N}\)~Evaluation};
  \node (s_j)   [process, minimum width=4cm, below of=s_y,  yshift=.5cm] {\(f, \bm{c}, \Jtot\)~Evaluation};
  \node (s_doe) [process, minimum width=4cm, below of=s_j,  yshift=.5cm] {Update~\(\overline{\D}_{J_i}, \overline{\D}_\text{sys}\)};
  \draw [arrow] (ps)  -- (s_y);
  \draw [arrow] (s_y) -- (s_j);
  \draw [arrow] (s_j) -- (s_doe);

  \node (exit) [
    decision, below of=s_doe, yshift=-1cm,
    minimum width=2.5cm, minimum height=.75cm, shape aspect=2,
  ] {Stop?};

  \node (outp) [io, below of=exit, yshift=0.25cm] {Optimization Output};

  \draw [arrow] (s_doe) -- (exit);

  \draw [arrow] (ss1_doe1.south) -- ([yshift=-0.5cm] ss1_doe1.south)
  -| ($ (inN_b.east)!0.5!(ins.west) $)
  -- (ins.west);

  \draw [arrow, red]
  (exit.west) node[anchor=south east] {no} -| ([xshift=-0.5cm]ss1_doe1.south west) |- (in1);

  \begin{pgfonlayer}{back1}
    \draw [arrow] (doe2_b.south) -- ([yshift=-0.5cm] doe2_b.south)
    -| node[anchor=south, rotate=90, xshift=\layerdy] {\(\bm{\cdots}\)} ($ (inN_b.east)!0.5!(ins.west) $)
    -- (ins.west);

    \draw [arrow, red]
    (exit.west) -| ([xshift=-0.5cm]ss1_doe1.south west) |- (in2_b);
  \end{pgfonlayer}

  \begin{pgfonlayer}{back3}
    \draw [arrow] (doeN_b.south) -- ([yshift=-0.5cm] doeN_b.south)
    -| ($ (inN_b.east)!0.5!(ins.west) $)
    -- (ins.west);

    \draw [arrow, red]
    (exit.west) -| ([xshift=-0.5cm]ss1_doe1.south west) |- node[anchor=south, rotate=-90, xshift=\layerdy] {\(\bm{\cdots}\)} (inN_b);
  \end{pgfonlayer}

  \draw [arrow] (exit) -- node[anchor=east] {yes} (outp);

\end{tikzpicture}
  }%
  \caption{Process flow for the Bayesian Algorithm for Collaborative Optimization.}\label{fig:bacodiag}
\end{figure}
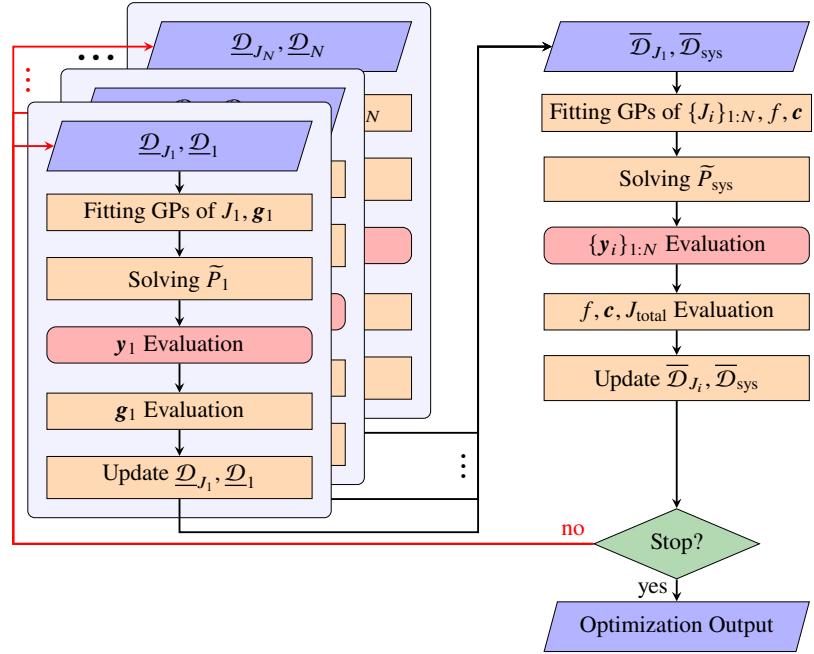

The two datasets share the same schema and are pooled when fitting the GP surrogate for~\(J_i\). Maintaining them separately preserves the indexing convention. Within a given iteration~\(k\), the subsystem level generates one new point appended to~\(\underline{\D_{J_i}}\), and the system level generates one new point appended to~\(\overline{\D_{J_i}}\). The combined size of~\(\overline{\D_{J_i}}^{(k)}\) and~\(\underline{\D_{J_i}}^{(k)}\) therefore grows by exactly two per iteration, and the superscript~\((k)\) on each dataset reflects the number of system-level iterations completed rather than the total number of evaluations. The subsystem constraint data are recorded in the following manner.
\begin{equation*}
  \D_i = \Bigl\{
    \overline{\bm{Y}}_{j\neq i},
    \underline{\bm{Z}}_i,
    \underline{\bm{X}}_i,
    \bm{g}_i\left(
      \overline{\bm{Y}}_{j\neq i},
      \underline{\bm{Z}}_i,
      \underline{\bm{X}}_i
    \right)
  \Bigr\},
\end{equation*}
where~\(\D_i\) records inputs with respect to~\(\ysys_{j \neq i}\) rather than the full~\(\ysys\), consistent with the argument structure of~\(\bm{g}_i\) in problem~\eqref{prob:cosubsys}.
The structural integration of the learning and optimization phases is visualized in \Cref{fig:bacodiag}. This diagram highlights the decoupling of the expensive function evaluations from the internal optimization loops.

\Cref{alg:baco} details BACO, depicting the enrichment of the GP surrogates at both the system and subsystem levels, the optimization of the acquisition functions, and the evaluation of the black-box functions at the proposed points.

\begin{algorithm}[!h]
  \caption{Bayesian Algorithm for Collaborative Optimization (BACO) Framework}\label{alg:baco}  \DontPrintSemicolon
  \SetKwInOut{Input}{Inputs}\SetKwInOut{Output}{Outputs}

  \Input{
    \ul{Functions}: \(f\): Objective function, \(\bm{c}\): Vector of system constraints, \(\bm{g}_i\): Vector of subsystem \(i\) constraints, \(J_i\): Subsystem \(i\) consistency discrepancy (for each \(i\))\\

    \ul{Initial DoEs}: \(\D_\text{sys}^{(0)}\), \(\D_i^{(0)}\), \(\overline{\D_{J_i}}^{(0)}\), \(\underline{\D_{J_i}}^{(0)}\) (for each \(i\))\\


  }
  \Output{\(\bm{z}_\text{best}, \bm{x}_\text{best}, \bm{y}_\text{best}, \fbst\).}

  \(k \gets 0\)

  \Repeat{
    a stopping criterion is met
  }{

    \ForEach{subsystem \(i \in \{ 1, \ldots, N \}\)}{

      \textcolor{green}{\tcp{--- Subsystem Level, Parallel Execution ---}}

      Fit GP surrogate for~\(J_i\) using \(
        \overline{\D_{J_i}}^{(k)}
        \cup
        \underline{\D_{J_i}}^{(k)}
      \), and for \(\bm{g}_i\) using \(\D_i^{(k)}\)

      Define acquisition function~\(\alpha_{J_i}\) using~\Cref{eqn:logei}

      Solve problem~\eqref{prob:bacosubsys}:\vspace{-.25cm}\textcolor{blue}{\[
          (
            \zsub_i^{(k+1)},
            \xsub_i^{(k+1)}
          )
          \in
          \argmax_{
            \zsub_i, \xsub_i
          }
          \left\{
            \alpha_{J_i}\left(
              \zsys^{(k)},
              \xsys_i^{(k)},
              \ysys^{(k)},
              \zsub_i,
              \xsub_i
            \right)
            :
            \bm{\mu}_{g_i}\left(
              \ysys_{j\neq i}^{(k)},
              \zsub_i,
              \xsub_i
            \right)\geq\bm{0}
          \right\}
      \]\vspace{-.5cm}}

      Evaluate and update DoEs:\vspace{-.25cm}
      \begin{alignat*}{3}
        \underline{\D_{J_i}}^{(k+1)} &\gets
        \underline{\D_{J_i}}^{(k)} &&\cup \{
          (
            \zsys^{(k)},
            \xsys_i^{(k)},
            \ysys^{(k)},
            \zsub_i^{(k+1)},
            \xsub_i^{(k+1)}
          )&&~,~
          J_i(
            \zsys^{(k)},
            \xsys_i^{(k)},
            \ysys^{(k)},
            \zsub_i^{(k+1)},
            \xsub_i^{(k+1)}
          )
        \}\\
        \D_i^{(k+1)} &\gets \D_i^{(k)} &&\cup \{
          (
            \ysys_{j\neq i}^{(k)},
            \zsub_i^{(k+1)},
            \xsub_i^{(k+1)}
          )&&~,~
          \bm{g}_i(
            \ysys_{j\neq i}^{(k)},
            \zsub_i^{(k+1)},
            \xsub_i^{(k+1)}
          )
        \}
      \end{alignat*}\vspace{-.5cm}
    }

    \textcolor{green}{\tcp{--- System Level ---}}

    Fit GP surrogates for~\(f\) and~\(\bm{c}\) using \(\D_\text{sys}^{(k)}\), and for~\(J_i\) using \(\overline{\D_{J_i}}^{(k)} \cup \underline{\D_{J_i}}^{(k+1)}\) (\(i=1:N\))

    Define acquisition function~\(\alpha_f\) using~\Cref{eqn:logei}

    Solve problem~\eqref{prob:bacosys}:\vspace{-.25cm}\textcolor{blue}{\[
        (
          \zsys^{(k+1)},
          \xsys^{(k+1)},
          \ysys^{(k+1)}
        )
        \in
        \argmax_{
          \zsys, \xsys, \ysys
        }
        \left\{
          \alpha_f(
            \zsys,
            \xsys,
            \ysys
          )
          :~
          \begin{matrix}
            \bm{\mu}_c(\zsys, \xsys, \ysys) \geq \bm{0} \\
            \mu_{J_i}(\zsys, \xsys_i, \ysys, \zsub_i^{(k+1)}, \xsub_i^{(k+1)}) = 0, \forall i
          \end{matrix}
        \right\}
    \]\vspace{-.5cm}}

    Evaluate and update DoEs:\vspace{-.25cm}
    \begin{alignat*}{6}
      &\D_\text{sys}^{(k+1)}&&\gets\D_\text{sys}^{(k)}&&\cup\{
        (
          \zsys^{(k+1)},
          \xsys^{(k+1)},
          \ysys^{(k+1)}
        ),~
        f(
          \zsys^{(k+1)},
          \xsys^{(k+1)},
          \ysys^{(k+1)}
        ),~\bm{c}(
          \zsys^{(k+1)},
          \xsys^{(k+1)},
          \ysys^{(k+1)}
        )
      \}\\
      &\overline{\D_{J_i}}^{(k+1)}&&\gets\overline{\D_{J_i}}^{(k)}&&\cup\{
        {(
            \zsys^{(k+1)},
            \xsys_i^{(k+1)},
            \ysys^{(k+1)},
            \zsub_i^{(k+1)},
            \xsub_i^{(k+1)}
        )},~J_i{(
            \zsys^{(k+1)},
            \xsys_i^{(k+1)},
            \ysys^{(k+1)},
            \zsub_i^{(k+1)},
            \xsub_i^{(k+1)}
        )}
      \}
    \end{alignat*}\vspace{-.5cm}

    Update best solution: \(\bm{z}_\text{best}, \bm{x}_\text{best}, \bm{y}_\text{best}, \fbst\).

    \(k \gets k + 1\)
  }

\end{algorithm}

\section{Black-box Evaluations and Stopping Criteria}\label{sec:conv}

\subsection{Black-box Evaluations}\label{subsec:eval}

In the context of the frameworks presented in this paper, one black-box evaluation refers to a single complete execution of a disciplinary simulation to map a specific set of inputs to the output vector \(\bm{y}_i\). It is designated as a black-box because the optimization algorithm only observes this input-output numerical relationship and cannot exploit the underlying algebraic structure or exact internal analytical derivatives of the governing physical equations. This is typically the most computationally expensive operation in the overall design process.

In CO, the computational cost is fundamentally dictated by the evaluation frequency of \(\bm{y}_i\). Because the discrepancy formulation mandates that \(\bm{y}_i\) is evaluated as a explicit function of the active system-level coupling variables, \(\ysys_{j\neq i}\), these expensive function calls cannot be isolated strictly to the lower-level subsystem optimizations. Every iteration of the upper system-level optimizer that perturbs the target states \(\ysys\) forces a simultaneous re-evaluation of \(\bm{y}_i\) across all \(N\) subsystems to accurately compute the discrepancy constraint \(J\). Letting \(K\) denote the total number of major iterations to reach convergence, \(L\) the average number of system-level optimization iterations per major iteration, and \(M_i\) the average number of local evaluations performed by the \(i\)-th subsystem optimizer per major iteration, the total computational complexity in terms of black-box calls scales as \(\Ord\left(K \left( N L + \sum_{i=1}^{N} M_i \right)\right)\). This strict coupling indicates that for expensive-to-compute multidisciplinary problems, the computational burden is heavily amplified by the dimensionality of the coupling variables, as the upper-level search directly multiplies the required number of high-fidelity simulation runs.

In BACO, nested optimization searches are replaced with GP surrogate evaluations to reduce computational cost. During each major iteration \(k\), the subsystem-level optimizers maximize an acquisition function \(\alpha_{J_i}\) subject to surrogate constraint predictions \(\bm{\mu}_{g_i}\). This internal search relies exclusively on surrogates, requiring no evaluations of the true disciplinary black-box functions \(\bm{y}_i\). Instead, \(\bm{y}_i\) is evaluated exactly once per subsystem after local optimization terminates to establish the true discrepancy and constraint values, which update the data sets \(\underline{\D_{J_i}}^{(k+1)}\) and \(\D_i^{(k+1)}\). Similarly, the system-level optimizer manipulates the coupling variables \(\ysys\) using the surrogate models \(\alpha_f\), \(\bm{\mu}_c\), and \(\mu_{J_i}\). Following this system-level update, \(\bm{y}_i\) is evaluated a second time across all \(N\) subsystems to determine the true discrepancy at the new system state, augmenting the upper-level data set \(\overline{\D_{J_i}}^{(k+1)}\). Because all iterative search procedures are performed on the GP surrogates, the framework strictly dictates two disciplinary evaluations per discipline per major iteration. For \(K\) major iterations, the exact total number of true black-box calls is \(2KN\), corresponding to an asymptotic computational complexity of \(\Ord(KN)\).

\subsection{Discrepancy and Constraint Violation Metrics}\label{subsec:metrics}

The state of a CO iteration is characterized by two scalar metrics evaluated at the system level after each subsystem optimization cycle. For a given system-level iterate~\((\zsys, \xsys, \ysys, \zsub^*, \xsub^*)\), the total discrepancy \(\Jtot\) aggregates the inter-disciplinary incompatibility across all subsystems,
\begin{equation}
  \Jtot\left(
    \zsys,
    \xsys,
    \ysys,
    \zsub^*,
    \xsub^*
  \right) = \sum_{i=1}^N J\left(
    \zsys,
    \xsys_i,
    \ysys,
    \zsub_i^*,
    \xsub_i^*
  \right),\label{eqn:disc_total}
\end{equation}
while the total constraint violation~\cite{audet2017}, \(\htot\), aggregates both subsystem-level and system-level infeasibility,
\begin{equation}
  \htot\left(
    \zsys,
    \xsys,
    \ysys,
    \zsub^*,
    \xsub^*
  \right)
  = \sum_{i=1}^N \sum_{l=1}^{m_i} \max\left\{
    0, - g_{i_l}\left(
      \ysys_{j\neq i},
      \zsub_i^*,
      \xsub_i^*
    \right)
  \right\}^{2}
  + \sum_{l=1}^m \max\left\{
    0, - c_l\left(
      \zsys,
      \xsys,
      \ysys
    \right)
  \right\}^{2}.\label{eqn:constviol_total}
\end{equation}
Both metrics are non-negative by construction, with \(\Jtot = 0\) attained only at a fully compatible point and \(\htot = 0\) only at a feasible point.

These metrics serve a dual role. They provide a stopping criterion. The iteration terminates at step \(k\) when \(\Jtot(\zsys^{(k)}, \ldots) \leq \epsilon_J\) and \(\htot(\zsys^{(k)}, \ldots) \leq \epsilon_h\) for prescribed tolerances \(\epsilon_J, \epsilon_h > 0\). They also serve as the primary axes along which solution quality is assessed, complementing the objective value \(f\).

\section{Numerical Results}
\subsection{Comparison Solvers and Implementation Details }\label{sec:implementation}

All frameworks were implemented in Python\footnote{The library built for this work is available at \href{https://moebehfn.github.io/mdotoolbox}{\textcolor{blue}{\texttt{moebehfn.github.io/mdotoolbox}}}.}. The primary metric for comparison is the number of true black-box discipline evaluations rather than wall-clock time, as the benchmark functions serve as analytic surrogates for genuine black-box systems. Unless stated otherwise, the parameters described below are common to all benchmark problems.

To assert a solution that is both as feasible as possible and as compatible as possible, all frameworks progression enforced a sufficient decrease rule on $\htot$ then $\Jtot$ before assessing the improvement of the objective value $f$. This can be seen in the below plots where in some cases, the $f$ or $\Jtot$ of the accepted iterant increases allowing for a solution that yields a better feasibility.

The three bi-level architectures share a common nested structure in which a system-level coordinator sets consistency targets for shared and coupling variables while subsystem-level optimizers satisfy local disciplinary constraints. Both levels employ COBYQA~\cite{ragonneau2023} as the local solver. Fifty percent of the total evaluation budget is reserved for the system level and the remainder is divided equally among the disciplines; an iteration ratio of 0.05 restricts each outer iteration to at most 5\% of the total budget across all levels. The total budget is 300 black-box evaluations. Convergence is declared when both~\(\Jtot\) and~\(\htot\) fall below~\(\varepsilon = 10^{-3}\). All four variants minimize the same subsystem-level weighted sum-of-squares discrepancy~\(J_i\) and differ only in how system-level coupling consistency is enforced: CO imposes~\(J_i = 0\); MCO replaces the shared design variables with the arithmetic mean of the subsystem copies so that the coordinator optimizes only the local and coupling targets; ICO augments the system-level objective with the penalty~\(\gamma \sum_i |J_i - r^2|\) using~\(\gamma = 1\),~\(r = 0\), and penalty update factor~\(\delta = 1.1\).

BACO replaces the direct optimizer calls at both levels with an acquisition function maximization loop built on GP surrogates constructed using the Surrogate Modeling Toolbox (SMT)~\cite{saves2024}. All GP models use KPLSK~\cite{bouhlel2019}, which applies partial least squares dimensionality reduction to the input space before fitting a stationary squared-exponential kernel in the reduced space. Length-scale hyperparameters are calibrated by maximizing the concentrated log-likelihood via COBYLA~\cite{powell1994} with 10 random restarts. At each iteration, the next evaluation point for each subsystem is selected by maximizing the log expected improvement acquisition function subject to GP-predicted subsystem feasibility. This maximization is performed using 25 multi-starts: 24 points drawn by LHS and 1 starting from the previous best solution. The same procedure is applied at the system level, with the additional surrogate-predicted constraint~\(\mu_{J_i} \leq \varepsilon\) on each discipline. BACO operates under the same overall evaluation budgets; all surrogate fitting and acquisition optimization are internal and do not count against the budget. Termination occurs when the budget is exhausted or both~\(\Jtot \leq \varepsilon\) and~\(\htot \leq \varepsilon\) are satisfied simultaneously with~\(\varepsilon = 10^{-3}\).

\subsection{Benchmarking on a Scalable MDO Problems}\label{sec:bench-results}

All convergence histories are plotted as a median trace across runs with a shaded 95\% confidence interval. The subplots report the objective value $f$, the total constraint violation $\htot$, and the total discrepancy $\Jtot$.

BACO is compared against CO, MCO, and ICO on a two-discipline instance of the Scalable MDO problem~\cite{tedford2010}, described in~\Cref{app:scalable}.

The stochastic nature of BACO, arising from the initial LHS DoE and the multi-start maximization of the acquisition function, as well as the sensitivity of gradient-free optimizers to initial starting points, necessitates a statistical treatment of performance. Each solver is launched from ten independent random starting points within the design space.

\begin{figure}[!ht]
  \centering
  \includegraphics[width=\textwidth]{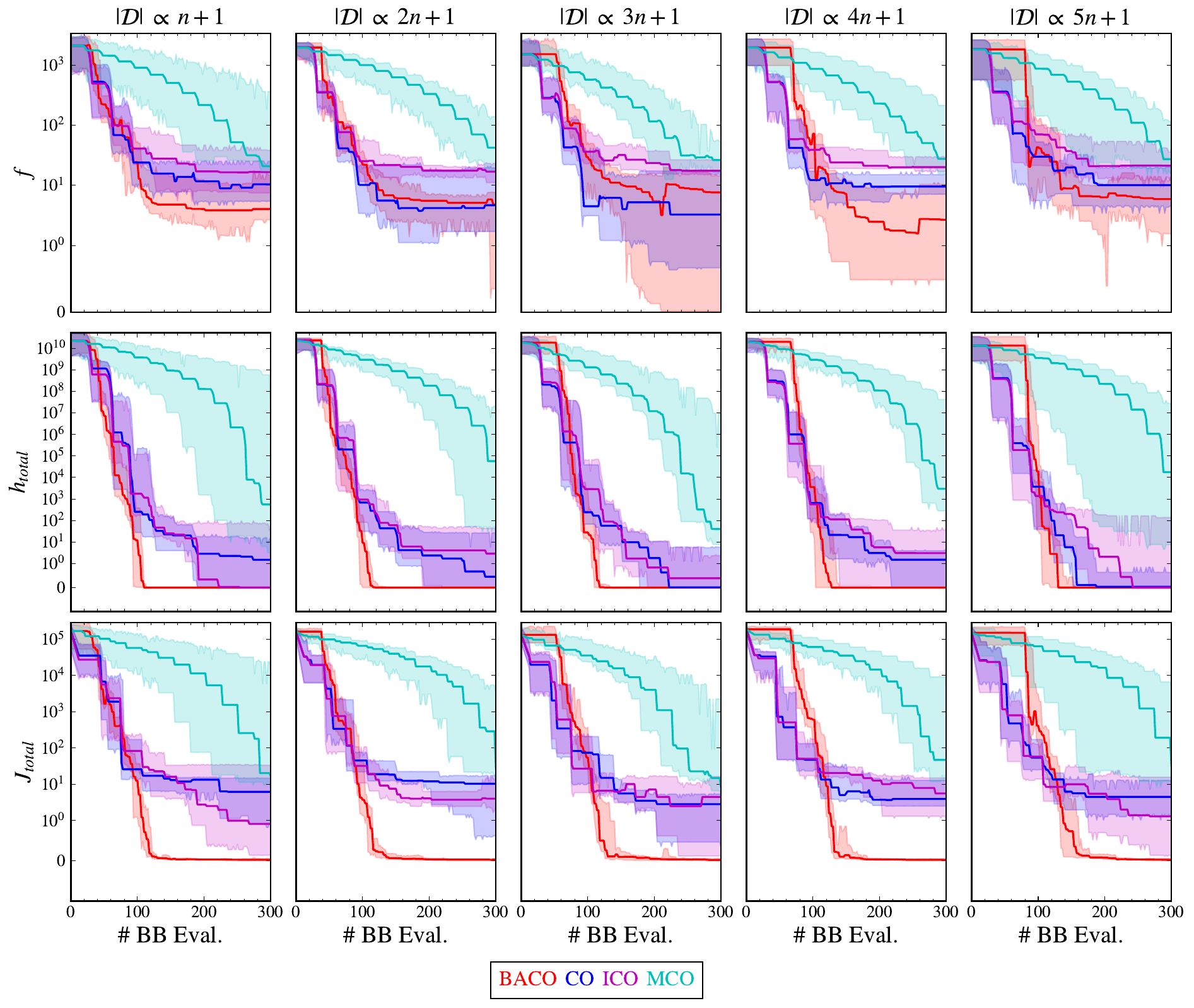}
  \caption{Aggregated convergence histories for the Scalable problem (10 runs per DoE size). DoE sizes are a function of $n$, the size of the system and the subsystem problems, depending on the targeted problem. On the x-axis label, \# BB Eval. indicates the number of black-box evaluations.}\label{fig:scalable}
\end{figure}

To assess the sensitivity of BACO to the initial DoE size, five DoE levels proportional to the problem dimension $n$ are evaluated: $n+1$, $2n+1$, $3n+1$, $4n+1$, and $5n+1$, where $n$ is the number of design variables of the targeted problem (system or subsystems). A separate LHS DoE is drawn for each subsystem and for the system coordinator at each level. The ten-run aggregated convergence results are summarized in~\Cref{fig:scalable}.

Across all DoE sizes, BACO achieves the lowest median objective value and drives both $\htot$ and $\Jtot$ to near-zero within the evaluation budget, outperforming all three CO variants. MCO exhibits the slowest convergence and the largest inter-run variance on all three metrics. CO and ICO perform comparably and occupy an intermediate range. The primary effect of increasing the DoE size is to shift BACO's onset of convergence rightward, as more evaluations are consumed during initialization, but the asymptotic solution quality is maintained or marginally improved due to a better-conditioned initial surrogate. These results indicate that BACO is robust to the choice of DoE size within the range tested, and that a DoE of $n+1$ points is sufficient to achieve competitive performance on this problem. Larger DoE sizes may be beneficial for more complex problems or those with higher-dimensional design spaces, but the marginal gains should be weighed against the increased initial evaluation cost.

\subsection{Results on a Coupled Aerodynamic Structural Problem}\label{sec:aeros-results}
The effectiveness of the BACO framework is further evaluated on a realistic aero-structural MDO problem. The problem involves the optimization of a generic transport aircraft wing, based on the CRM geometry, to minimize fuel burn while satisfying aerodynamic and structural constraints. The system-level objective is to minimize the total fuel burn during a cruise mission, calculated using the Breguet range equation~\cite{coffin1920,raymer2018,anderson1999}. The aircraft operates at Mach \num{0.84} at \SI{10000}{\meter} altitude. The total weight is the sum of the fixed operating empty weight, payload weight, and structural wing weight. The target range is~$R =$~\SI{10000}{\kilo\meter} and the thrust-specific fuel consumption is~$c_T =$~\SI{1.5e-5}{\kilo\gram\per\newton\per\second}. Fuel burn is computed as~\cite{coffin1920,raymer2018}:
\begin{equation*}
  W_f = (W_0 + W_s)\left(\exp\!\left(\frac{R c_T g}{V_\infty C_L/C_D}\right) - 1\right).
\end{equation*}

The problem is decomposed into aerodynamics and structures\footnote{The aerodynamic and structural simulations were conducted using a custom-built solver available at \href{https://moebehfn.github.io/lightaero}{\textcolor{blue}{\texttt{moebehfn.github.io/lightaero}}}.} coupled disciplines.

The aerodynamic discipline uses a Vortex Lattice Method (VLM)~\cite{falkner1943} solver to compute the lift coefficient $C_L$, drag coefficient $C_D$, and the span-wise sectional lift distribution $c_l(y)$. The wing is discretized into 10 span-wise panels. The effective local angle of attack accounts for geometric twist $\theta(y)$ and an aero-elastic washout correction proportional to the tip deflection ratio $\delta/b$ supplied by the structural discipline. The structural discipline models the wing spar as a hollow circular Euler-Bernoulli beam~\cite{hodges2011}, discretized at the same span-wise stations as the aerodynamic panels. Sectional lift coefficients from the aerodynamic discipline are converted to distributed aerodynamic loads, from which tip deflection and internal bending stresses are recovered via a finite-element beam solve. The problem treats half of the wing with 5 span-wise panels, and the results are mirrored to the other half. The shared design variable is $z = \alpha$, the local variables $x_1$ and $x_2$ are the span-wise twist and spar thickness distributions, and the coupling variables are $y_1 = (C_L, C_D, c_l(y))$ and $y_2 = (\delta/b, W_s)$. All variables are normalized to~$[0,1]$ using their physical bounds shown in \Cref{tab:aerostruct_bounds}. The problem is formulated as follows:
\begin{mini*}
  {\zsys, \xsys, \ysys}{W_f}{}{}
  \addConstraint{C_L q_\infty S_\text{ref}}{\geq W_0 + W_s}{\quad \text{(lift)}}
  \addConstraint{\delta/b}{\leq \num{0.1}}{\quad \text{(deflection)}}
  \addConstraint{C_L}{\geq \num{0.1}}{\quad \text{(stall)}}
  \addConstraint{\num{2.5}\times\sigma_j}{\leq \sigma_Y,~j=1,\ldots,5}{\quad \text{(stress)}}
  \addConstraint{\bm{y}_1}{= \mathcal{A}(\bm{z}, x_1, y_2)}{\quad \text{(aerodynamics)}}
  \addConstraint{y_2}{= \mathcal{S}(\bm{z}, x_2, y_1)}{\quad \text{(structures)}}
\end{mini*}
\begin{table}[!h]
  \centering
  \caption{Physical bounds for design and coupling variables.}
  \label{tab:aerostruct_bounds}
  \begin{tabular}{ll@{~}llrrr}
    \hline
    Variable & & & Description & Lower & Upper & Initial Guess \\
    \hline
    $\alpha$ & $z$ & $\in \R$ & Angle of attack & \ang{-2} & \ang{+8} & \ang{3} \\
    $\theta$ & $\bm{x}_1$ & $\in \R^5$ & Twist & \ang{-10} & \ang{+10} & \{\ang{-1}\ldots\ang{+1}\ldots\ang{-1}\} \\
    $t$ & $\bm{x}_2$ & $\in \R^5$ & Spar thickness & \SI{5}{\milli\meter} & \SI{100}{\milli\meter} & \SI{20}{\milli\meter} everywhere \\
    $C_L$ & $[\bm{y}_1]_1$ & $\in \R$ & Lift coefficient & \num{0.05} & \num{1.20} & \num{0.6} \\
    $C_D$ & $[\bm{y}_1]_2$ & $\in \R$ & Drag coefficient & \num{0.012} & \num{0.100} & \num{0.05} \\
    $c_{l,i}$ & $[\bm{y}_1]_{3:7}$ & $\in \R^5$ & Sectional lift coefficient & \num{0.0} & \num{3.0} & \num{0.6} everywhere \\
    $\delta/b$ & $[\bm{y}_2]_1$ & $\in \R$ & Tip deflection ratio & \num{0.0} & \num{0.2} & \num{0.05} \\
    $W_s$ & $[\bm{y}_2]_2$ & $\in \R$ & Structural weight & \SI{1e4}{\newton} & \SI{6e5}{\newton} & \SI{1e5}{\newton} \\
    \hline
  \end{tabular}
\end{table}

\begin{figure}[!h]
  \centering
  \includegraphics[width=0.8\textwidth]{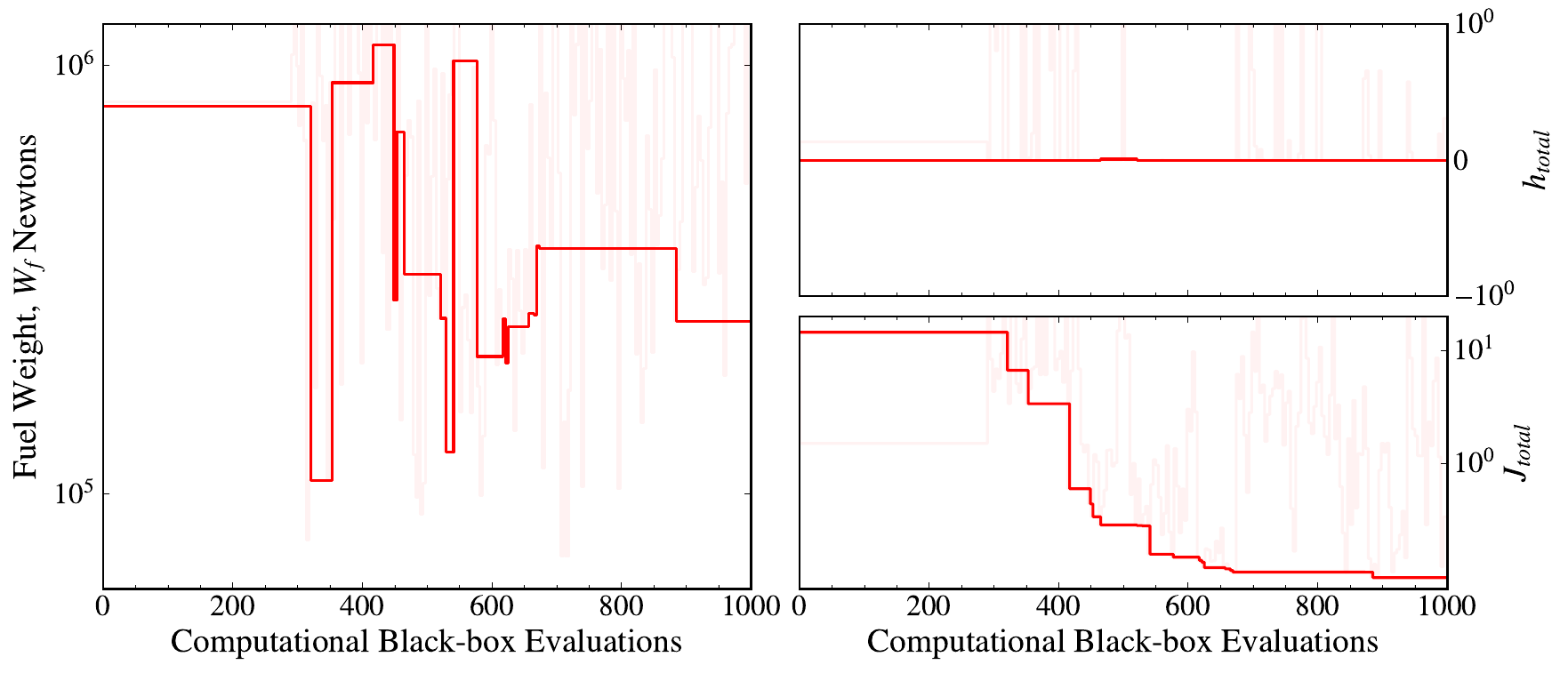}
  \caption{Convergence history for the aero-structural optimization problem.}\label{fig:ucrm_history}
\end{figure}

BACO was allocated a budget of 1000 evaluations for this problem. Each subsystem surrogate is seeded by LHS with a DoE size of~\(5n + 1\) following the sensitivity analysis done in \Cref{sec:bench-results}. The convergence history is depicted in \Cref{fig:ucrm_history}.

BACO accepted 18 improvement over the allocated budget, with the final accepted solution at eval 886 with full feasibility and a residual coupling discrepancy of $\Jtot = 0.089$. Feasibility is established early and largely maintained thereafter, with $\htot$ remaining at or near zero for the majority of system evaluations. The objective value of accepted iterants is non-monotone due to the hierarchical ordering imposed by BACO, later accepted points may carry a worse $f$ than earlier ones if they are accepted on the basis of improved interdisciplinary compatibility alone. $\Jtot$ decreases monotonically across accepted iterants, reflecting consistent tightening of the coupling residual toward disciplinary equilibrium. The run stagnates towards the final remaining evaluations.

A comparison of the initial and optimized solutions is provided in \Cref{tab:ucrm}.

\begin{table}[!h]
  \centering
  \caption{Initial and optimized solutions for the aero-structural optimization problem.}\label{tab:ucrm}
  \begin{tabular}{llllllll}
    \toprule
    & Fuel  & Total & Total & Lift & Drag & Glide & Angle of \\
    & Weight & Constraint & Discrepancy & Coefficient & Coefficient & Ratio & Attack \\
    & [\si{\kilo\newton}] & Violation & & & [Counts] & & [\si{\degree}] \\
    \midrule
    Initial   & \num{803.0} & 0 & 14.517 & 0.324 & 97.167 & 33.292 & 3.0 \\
    Optimized & \num{252.25} & 0 & 0.089 & 0.512 & 201.685 & 25.370 & 0.105 \\
    \bottomrule
  \end{tabular}
\end{table}

The optimized solution achieves a significant reduction in fuel burn, from  \SI{803.0}{\kilo\newton} to \SI{252.25}{\kilo\newton}, while reducing the constraint violation by three orders of magnitude. The lift and drag coefficients both increase bringing the glide ratio, $L/D$, to a more typical value for a transport aircraft of similar wing dimensions.

\begin{figure}[!h]
  \centering
  \begin{subfigure}[t]{.48\textwidth}
    \includegraphics[width=\textwidth]{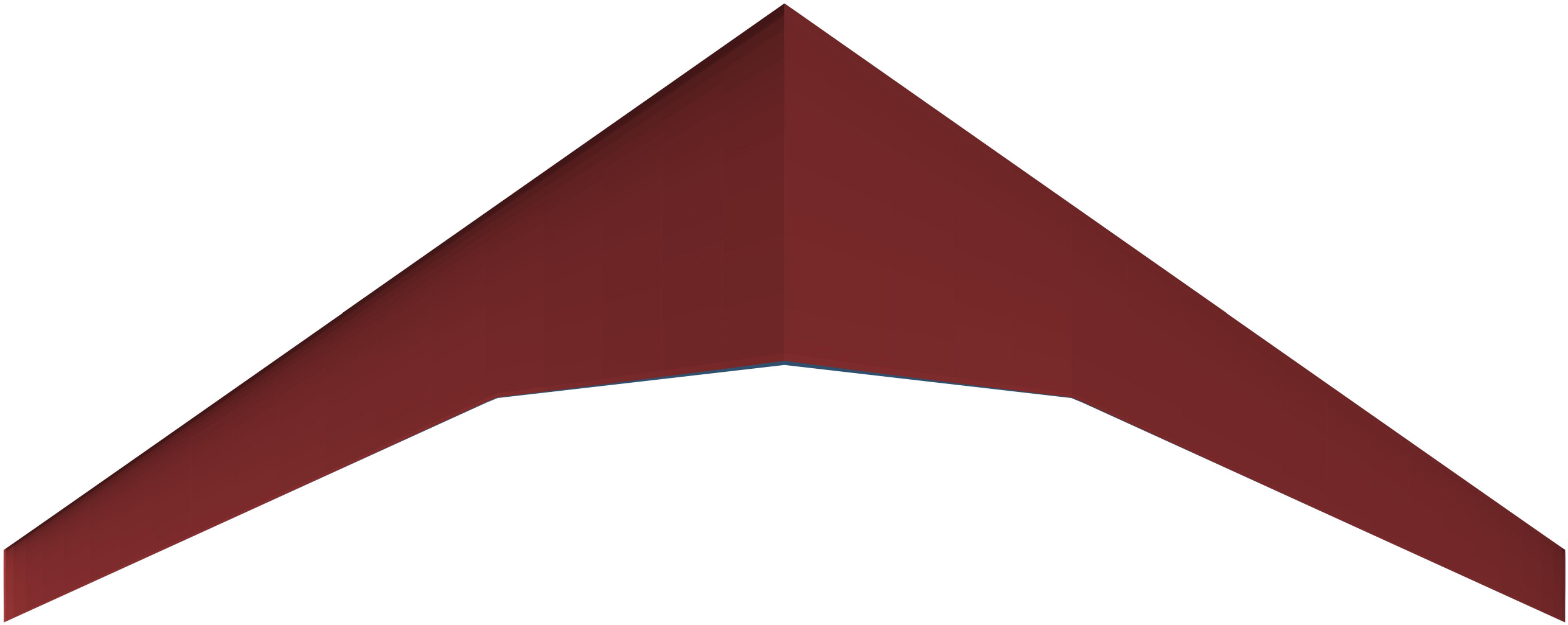}
    \caption{Top view}\label{fig:wing-top}
  \end{subfigure}
  \hfill
  \begin{subfigure}[t]{.48\textwidth}
    \includegraphics[width=\textwidth]{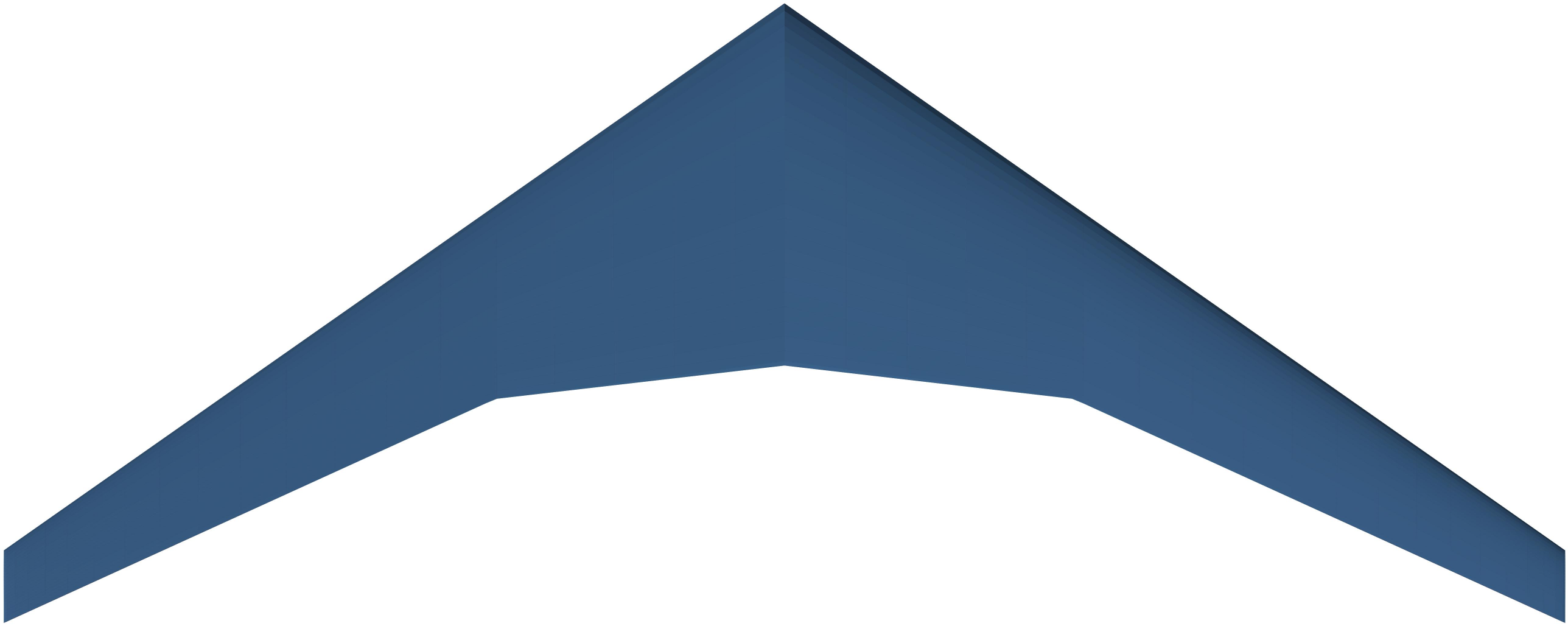}
    \caption{Bottom view}\label{fig:wing-bottom}
  \end{subfigure}
  \vspace{.25cm}
  \begin{subfigure}[t]{.48\textwidth}
    \includegraphics[width=\textwidth]{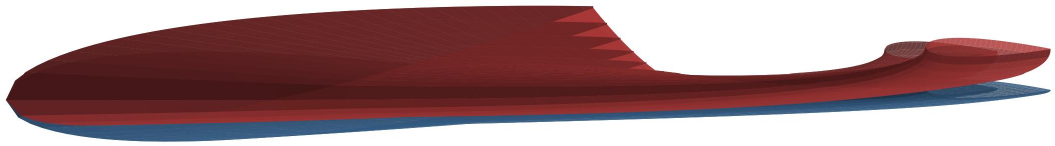}
    \caption{Left view}\label{fig:wing-left}
  \end{subfigure}
  \hfill
  \begin{subfigure}[t]{.48\textwidth}
    \includegraphics[width=\textwidth]{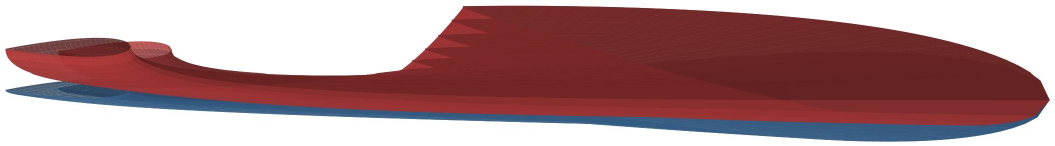}
    \caption{Right view}\label{fig:wing-right}
  \end{subfigure}
  \vspace{.25cm}
  \begin{subfigure}[t]{.48\textwidth}
    \includegraphics[width=\textwidth]{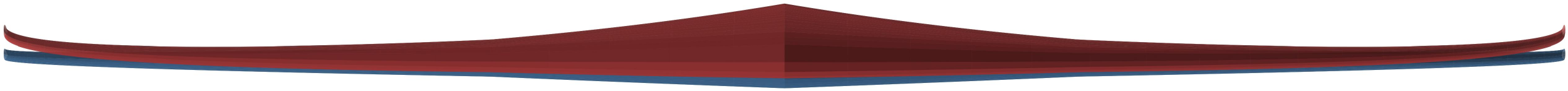}
    \caption{Front view}\label{fig:wing-front}
  \end{subfigure}
  \hfill
  \begin{subfigure}[t]{.48\textwidth}
    \includegraphics[width=\textwidth]{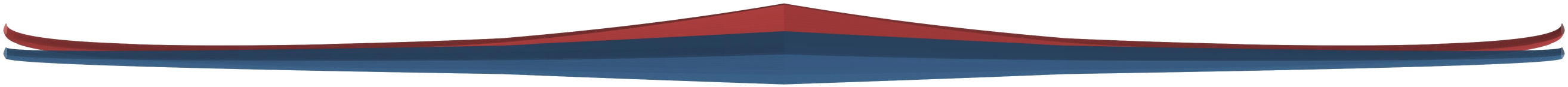}
    \caption{Rear view}\label{fig:wing-rear}
  \end{subfigure}
  \caption{Initial (blue) and optimized (red) wing. Wing tip deflection applied as a cubic interpolation of the wing tip deflection.}\label{fig:wing}
\end{figure}

The produced wing overlain on the initial wing can be seen in \Cref{fig:wing}.

The top and bottom views in \Cref{fig:wing-top,fig:wing-bottom} show separation between the two geometries, with the top view dominated by the optimized wing and the bottom view by the initial configuration. This is indicative of the shift in the twist distribution across the span.

The side views in \Cref{fig:wing-left,fig:wing-right} emphasize the difference in span-wise twist, with the optimized wing exhibiting a more pronounced nose-down washout near the root and a flatter profile toward the tip, consistent with the redistribution of lift away from the inboard stations to reduce root bending stress. The cubic interpolation at the tip produces a smooth geometric transition that is clearly distinguishable from the initial configuration. The wingtip deflection due to the span-loads is also visible in both views on the optimized wing.

The front and rear views in \Cref{fig:wing-front,fig:wing-rear} indicate that the optimized wing sustains greater tip deflection relative to the initial configuration, suggesting the structural solution operates closer to the tip deflection constraint $\delta/b \leq 0.1$, which is consistent with a reduced spar thickness distribution that lowers structural weight $W_s$ and improves fuel burn.

\begin{figure}[!h]
  \centering
  \includegraphics[width=.8\textwidth]{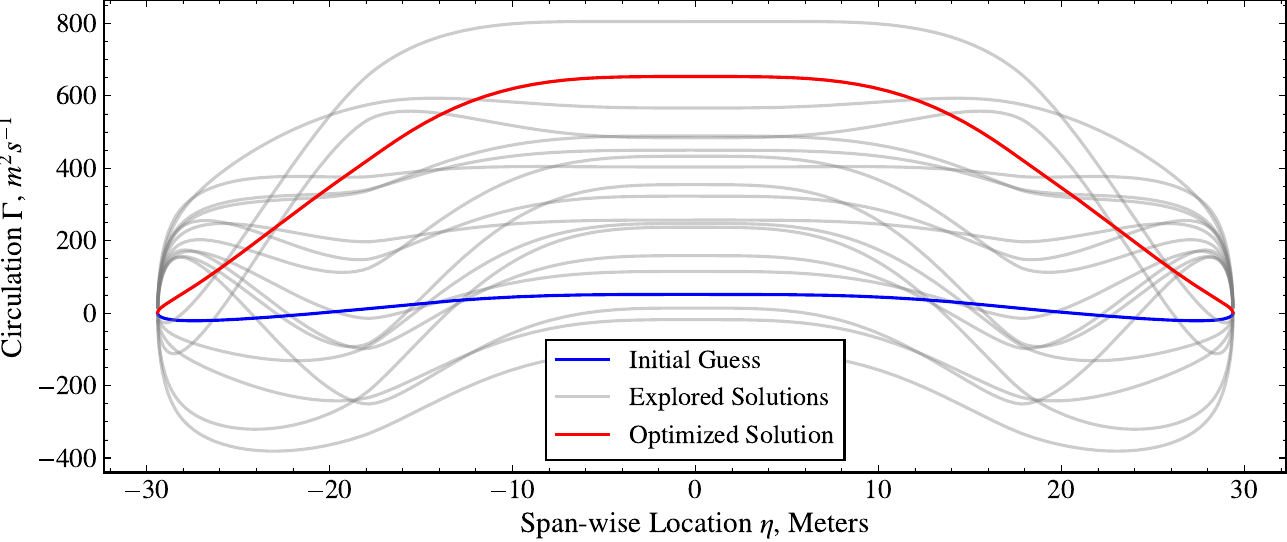}
  \caption{Span-wise circulation distribution $\Gamma(\eta)$ at the initial guess and the optimized solution. Grey traces depict the solutions marked as improving by BACO.}\label{fig:circulation}
\end{figure}

The span-wise circulation distribution across accepted iterants is shown in~\Cref{fig:circulation}.

The initial guess yields a near-zero circulation that is physically inconsistent with the required lift generation. The intermediate solutions span a wide range, including configurations with negative circulation over portions of the span, reflecting the exploratory behavior of BACO prior to surrogate refinement. The optimized circulation is smooth, symmetric, and approximately elliptic. The distribution peaks inboard of the mid-semi-span and exhibits a slightly flat-topped profile rather than a true semi-ellipse, which is consistent with the active bending stress constraints which tends to redistribute lift away from the tip to reduce the root bending moment, allowing a reduction in spar thickness and structural weight at the cost of a small increase in induced drag.
\section{Conclusion}\label{sec:conclusion}

This paper presented BACO, a Bayesian optimization framework for CO that replaces nested black-box calls at both the system and subsystem levels with GP surrogates and acquisition function maximization. The key structural contribution is the dual-dataset formulation for each discrepancy function $J_i$, populated separately from subsystem-level and system-level evaluations and pooled for surrogate fitting. This allows the GP model of $J_i$ to receive exactly two true evaluations per major iteration regardless of problem dimensionality, reducing drastically the number of calls to the black-box .

Benchmarking on the Scalable MDO problem over 50 randomized instances across five DoE sizes demonstrated that BACO consistently achieves lower objective values and drives both the total constraint violation and the total discrepancy to near-zero within a competitive budget, outperforming the state-of-the-art and independently of the initial DoE size. On the aero-structural wing twist optimization problem, BACO identified a feasible solution with $\Jtot = 0.089$ within 886 of 1000 allocated evaluations, recovering a span-wise lift distribution that deviates from the elliptical optimum in a manner physically consistent with active bending stress and tip deflection constraints.

Future work will explore the application of BACO to a broader range of MDO problems, including those with higher-dimensional design spaces and more complex coupling structures~\cite{saves2024}. A prospective avenue is the extension of BACO to handle multi-fidelity settings~\cite{cordelier2026}. This is expected to further reduce the computational cost of MDO process.

\section{Scalable MDO Problem}\label{app:scalable}
The scalable MDO problem~\cite{tedford2010} was designed to test the effect of increasing dimensionality while maintaining manageable computational costs. It allows for varying the number of disciplines, coupling variables, local design variables, and global design variables.

For this study, a quadratic objective was used. The disciplines have a linear dependence on each other and feature a local constraint on each coupling variable. The governing equations for each discipline are a linear system dependent on global design variables, local design variables, and non-local coupling variables.

The optimization problem statement is~\cite{tedford2010}
\begin{align*}
  \min_{\bm{z},\bm{x}} \quad & \bm{z}^\top\bm{z}+\sum_{i=1}^{N}\bm{y}_i^\top\bm{y}_i \\
  \text{subject to:} \quad & \bm{1}-\bm{C}_i^{-1}\bm{y}_i\le\bm{0}, \quad i=1,\ldots,N
\end{align*}

The governing equations for each discipline~\(i\) are given by
\begin{equation*}
  \bm{y}_i(\bm{z},\bm{x}_i,\bm{y}_j)=-\bm{C}_i^{-1}(\bm{C}_{z}\bm{z}+\bm{C}_{x_i}\bm{x}_i-\bm{C}_j\bm{y}_j)
\end{equation*}
In this formulation, all~\(\bm{C}\) matrices are composed of random positive coefficients generated prior to the optimization. For this exercise,~\(\bm{C}_i\) matrices are unitary, following the example shown in~\cite[Section 5.4]{tedford2010}. The problem is formulated as a simple two-discipline optimization problem. The variable dimensions and bounds are summarized in the table below. The coefficient matrices \(\bm{C}\in\R^{1\times1}\)  are drawn from \([0, 10]\) with a fixed seed of 42.
\begin{table}[!h]
  \centering
  \begin{tabular}{lcccccc}
    \toprule
    &\(\bm{x}\) & \(x_1\) & \(x_2\) & \(y_1\) & \(y_2\) & Total \\
    \midrule
    Dimensions & 1 & 1 & 1 & 1 & 1 & 5 \\
    Bounds & \multicolumn{3}{c}{\([-10, +10]\)} & \multicolumn{2}{c}{\([-100, +100]\)} & -- \\
    \bottomrule
  \end{tabular}
\end{table}

\bibliography{bib}
\end{document}